\documentclass[12pt]{article}
\usepackage{amsthm, amsmath, amssymb,latexsym}
\newtheorem{theorem}{Theorem}[section]
\newtheorem{lemma}[theorem]{Lemma}
\newtheorem{corollary}[theorem]{Corollary}
\newtheorem{proposition}[theorem]{Proposition}
\newtheorem{conjecture}[theorem]{Conjecture}
\newtheorem{remark}[theorem]{Remark}
\newtheorem{definition}[theorem]{Definition}

\def\A{\mathcal{A}}
\def\O{\mathcal{O}}

\newcommand{\nc}{\newcommand}
\nc{\cH}{{\mathcal H}}
\nc{\cA}{{\mathcal A}}
\nc{\cG}{{\mathcal G}}
\nc{\cC}{{\mathcal C}}
\nc{\cO}{{\mathcal O}}
\nc{\cI}{{\mathcal I}}
\nc{\cB}{{\mathcal B}}
\nc{\cY}{{\mathcal Y}}
\nc{\cK}{{\mathcal K}}
\nc{\cX}{{\mathcal X}}
\nc{\cS}{{\mathcal S}}
\nc{\cE}{{\mathcal E}}
\nc{\cF}{{\mathcal F}}
\nc{\cZ}{{\mathcal Z}}
\nc{\cQ}{{\mathcal Q}}
\nc{\cN}{{\mathcal N}}
\nc{\cP}{{\mathcal P}}
\nc{\cL}{{\mathcal L}}
\nc{\cM}{{\mathcal M}}
\nc{\cT}{{\mathcal T}}
\nc{\cW}{{\mathcal W}}
\nc{\cU}{{\mathcal U}}
\nc{\cD}{{\mathcal D}}
\nc{\cJ}{{\mathcal J}}
\nc{\cV}{{\mathcal V}}
\nc{\bH}{{\mathbb H}}
\nc{\bA}{{\mathbb A}}
\nc{\bG}{{\mathbb G}}
\nc{\bC}{{\mathbb C}}
\nc{\bO}{{\mathbb O}}
\nc{\bI}{{\mathbb I}}
\nc{\bB}{{\mathbb B}}
\nc{\bY}{{\mathbb Y}}
\nc{\bK}{{\mathbb K}}
\nc{\bX}{{\mathbb X}}
\nc{\bS}{{\mathbb S}}
\nc{\bE}{{\mathbb E}}
\nc{\bF}{{\mathbb F}}
\nc{\bZ}{{\mathbb Z}}
\nc{\bQ}{{\mathbb Q}}
\nc{\bN}{{\mathbb N}}
\nc{\bP}{{\mathbb P}}
\nc{\bL}{{\mathbb L}}
\nc{\bM}{{\mathbb M}}
\nc{\bT}{{\mathbb T}}
\nc{\bW}{{\mathbb W}}
\nc{\bU}{{\mathbb U}}
\nc{\bD}{{\mathbb D}}
\nc{\bJ}{{\mathbb J}}
\nc{\bV}{{\mathbb V}}
\nc{\bbZ}{{\mathbb Z}}
\nc{\bR}{{\mathbb R}}
\nc{\til}{\tilde{X}}
\nc{\fr}{{\rightarrow}}
\nc{\co}{{\overline{\nabla}}}
\nc{\lra}{\longrightarrow}
\nc{\ra}{\rightarrow}

\begin{document}

                                %
                                %
                                %

                                %
                                %
                          %
                                %
                                %
                                %
\title{Osculating spaces and diophantine equations\\
\small{with an Appendix by Pietro Corvaja and Umberto Zannier.}}

\author{ Michele Bolognesi, Gian Pietro Pirola}
\footnotetext {Partially supported by
 1) PRIN 2005  {\em "Spazi di moduli e teoria di Lie"}; 2) Indam Gnsaga;
3) Far 2006 (PV):
 {\em Variet\`{a} algebriche, calcolo
 algebrico,
 grafi orientati e topologici}.}
                                %
                                %
\maketitle
                                %
                                %
 \begin{abstract}
{\em } \vskip3mm

This paper deals with some classical problems about the projective geometry of complex algebraic curves.
We call \textit{locally toric} a projective curve that in a neighbourhood of every point has a local analytical
parametrization of type $(t^{a_1},\dots,t^{a_n})$, with $a_1,\dots, a_n$ relatively prime positive integers. 
In this paper we prove that the general tangent line to a locally toric curve in $\bP^3$ meets the curve only at the 
point of tangency. This result extends and simplifies those of the paper \cite{kaji} by H.Kaji where the same 
result is 
proven for any curve in $\bP^3$ such that every branch is smooth. More generally, under mild hypotesis, up to a 
finite number of anomalous parametrizations $(t^{a_1},\dots,t^{a_n})$, the general osculating 2-space to a locally 
toric curve of genus $g<2$ in $\bP^4$ does not meet the curve again. The arithmetic part of the proof of this result
relies on the Appendix \cite{cz:rk} to this paper. By means of the same methods we give some applications and 
we propose possible further developments.

\end{abstract}
                                %
                                %

\thispagestyle{empty}           %
\pagestyle{myheadings}
\markright{Osculating spaces and diophantine equations} %
                                %
                                %



                                %
                                %
                                %

\parindent=0pt

\section*{Introduction}

Let $C$ be a smooth connected complex complete curve and let

$$\phi: C\to \mathbb P^3$$

be a morphism such that $\phi$ is birational onto its image. In \cite{kaji} H.Kaji proved that if $\phi$ is
unramified then the tangent line to the general point does not meet
 again the curve. This gives a partial answer to a problem posed by  A.Terracini in a paper of 1932 \cite{terra}.
Our paper is a sort of revival of the beautiful Kaji's  argument, which is by way of contradiction. The proof of \cite{kaji} splits in two parts. In the first part one proves the existence of a flex point on $\phi(C)$ such that the tangent  line deforms to a trisecant
tangent line; the associated infinitesimal condition provides simpler equations. We noticed that this procedure
realizes a non-trivial reduction of the problem to the case of some special rational curves. This perspective allows us to generalize the viewpoint. In fact in \cite{kaji} the author
considered only curves $C\subset \bP^3$ such that the normalization morphism $C' \ra C$ is unramified, whereas we
allow also some ramification. More precisely, the curves that we will study are of the following type.

\begin{definition}\label{toric}
We say that $\phi:C\to \bP^n$ is locally toric if for any point of $C,$ there are n positive relatively prime
integers $a_i$, $i=1,...,n,$  such that $0<a_i<a_{i+1}$ and
a local analytical parametrization of $\phi$ of the form:

\begin{equation}  \label{affine}
(t^{a_1},...,t^{a_{n}}).
\end{equation}

\end{definition}

\begin{remark}
By relatively prime integers $a_i$ we always mean that $(a_1,...,a_n)=1$ as an ideal of $\bZ$.
\end{remark}

We will often speak of a \textit{locally toric curve} when referring to a curve with a locally toric morphism to a projective space $\bP^n$.
We remark that locally toric with $a_1=1$ means exactly that the normalization is unramified. Furthermore the
affine curve in (\ref{affine}) has a natural toric action. Our first main result is the following:

\begin{theorem}\label{principale}
Let $C$ be a connected smooth complete algebraic curve. If $\phi : C \to \bP^3$ is a morphism birational onto
its image that is locally toric then the tangent line to the general point does not meet again the curve.
\end{theorem}

The second part in \cite{kaji} is more computational and it relies on a classical, non-trivial result of Enestr\"{o}m and Kakeya
\cite{kake} about the zeros of a polynomial with real coefficients. We find a completely elementary proof of a
more general result as an application of Rolle's theorem. Furthermore we remark that in the case of locally toric curves the computation reduces to counting the integer zeroes of an exponential diophantine equation.

The paper is organized in the folllowing way.
In section \ref{uno} we prove our basic lemma. In section \ref{secred} we perform the Kaji reduction. As it was already pointed out by Kaji by giving explicit examples in positive characteristic, the hypothesis of
working in characteristic zero is essential. In fact if the general
tangent line were a trisecant then the principal parts bundle $\cP^1$ of $\cO_C(1)$ would split, which is impossible in
characteristic zero. The technical problem
in our case is that, while Kaji used the principal parts bundle, we need to build a vector bundle that is a
sub-sheaf of the principal parts bundle of $C$. In the third section we discuss some applications of the same
methods, studying the problem of whether the span of the tangent
lines at two general points of a locally toric curve in $\bP^4$ containes the tangent line to a third point. Then we discuss a more general
``dual problem" concerning the osculating linear spaces to a couple of general points. Finally in section 4 
we consider the problem of whether the $(n-2)$-linear osculating space at a general point of a locally toric
curve in $\bP^n$ intersects again the curve. If $g<2$, under some mild technical assumptions (see Thm. \ref{evev}) we are able to reduce this problem to that of calculating the rank of the following $n\times (n-1)$ matrix, where $z\in \bC$.

\begin{equation} \label{Amm}
A_{a_1,a_2,\dots,a_n}(z)=\begin{pmatrix}
 a_1  & a_1^2  &\dots & a_1^{n-2} & z^{a_1}-1 \\
a_2  & a_2^2  &\dots & a_2^{n-2}&z^{a_2}-1 \\
\dots & \dots   &\dots &\dots &\dots  \\
a_n  & a_n^2  &\dots & a_n^{n-2}&z^{a_n}-1
\end{pmatrix}
\end{equation}

In fact the general $(n-2)$-linear osculating space does not intersect again the curve if and only if the rank of matrix \ref{Amm} is $n-1$, for $z\neq 1 $. We will call \textit{anomalous} (see Definition \ref{ano} and \ref{anoc}) a $n$-plet of relatively prime integers $0<a_1< \dots < a_n$ such that $A_{a_1,\dots,a_n}(z)$ has rank smaller than $n-1$ for some $z\neq 1$ and \textit{anomalous locally toric curve} a locally toric curve with an anomalous $n$-plets of exponents in its parametrization.

\begin{remark}\label{cotte}
In the appendix to this paper \cite{cz:rk}, Pietro Corvaja and Umberto Zannier prove, via methods 
related to the paper \cite{bomaza}, that there
exists only a finite (possibly zero) number of anomalous $4$-plets $0<a_1< \dots < a_4$ of integers. In particular they find an explicit bound for $a_4$ thus reducing the problem of finding anomalous $4$-plets to a finite number of verifications. They believe that the same property of finiteness can be proved for every $n$.  Moreover they give an algorithm, that could be performed by a calculator, to explicitly compute all anomalous $4$-plets, but this has not been done yet. Performing this algorithm would answer the question of the existence of such $4$-plets, proving Conjecture $\ref{matrice}$ that we state in Section $\ref{uno}$.
\end{remark}

A (partial) answer to the problem of $(n-2)$-osculating spaces is our second main result. Let $\cP^n$ the bundle of principal parts
of order $n$, we have the following theorem.
 
\begin{theorem}\label{evev}
Let $C\subset \bP^n=\bP V$ be a locally toric curve and $g(C)=0$ or $1$. If $g(C)=1$, we suppose moreover that the natural evaluation
map

$$ev_{n-2}:V \longrightarrow \cP^{n-2}$$

is surjective. Then, if $C$ is not anomalous, the general 
osculating $(n-2)$-plane does not intersect the curve in a second point different from the osculating one.
\end{theorem}

We remark that, despite this generalization, the case of singular curves is still widely open.
We can formulate the following:

\begin{conjecture} Let $X\subset \bP^3$ be a non-degenerate singular complex curve.
Then the general tangent line does not meet again the curve.\end{conjecture}

Moreover, from the enumerative point of view, the problem of tangent lines intersecting again the curve $C\subset \bP^3$ has 
been studied in \cite{lb}, where the number of such lines is given, once one supposes that it is finite. We believe that section 1.2 of this paper could open interesting perspectives under an enumerative point of view as well.
 
The arithmetic part of section 4 relies heavily on the results by P.Corvaja and U.Zannier contained in the appendix \cite{cz:rk}. We would like to
thank Enrico Schlesinger for many fruitful conversations, Pietro Corvaja and Umberto Zannier for the intense and interesting
E-mail correspondence we had while writing this paper.

\textbf{Added in Proof:} We have been kindly informed by J.Starr that Theorem \ref{principale} has been proved indipendently with similar techniques by I.Coskun,
N.Elkies, G. Farkas, J.Harris and J.Starr. Moreover they show that Theorem \ref{principale} is no longer true if one considers analytic arcs.


\section{Ranks of polynomial matrices.}\label{uno}

\subsection{An elementary lemma}
 Let $a,b$ and $c$ be positive integers, with $0< a< b< c.$
 let $f=GCD(a,b,c)$, note that $f=1$ if and only if the numbers are relatively prime.
We take $z$ in the complex numbers field. Let us consider now the matrix:

\begin{equation} \label{Am}
A(z)=\begin{pmatrix}
 a  &\  z^{a}-1 \\
 b &\   z^{b}-1\\
 c &\  z^{c}-1 \\
\end{pmatrix}
\end{equation}

If we define $v(z):=(z^a, z^b, z^c)$, $u:=(1,1,1)$ and $w:=(a,b,c)$,
 we have that $A(z)$ has rank $1$ when $v-u$ is proportional to $w.$

\begin{lemma} \label{lui}
The rank of $A(z)$ is $1$ if and only if $z^f= 1.$
\end{lemma}
\begin{proof}
If the rank of $A(z)$ is $1$, both the real and the
immaginary part  of $(z^{a}-1,  z^{b}-1, z^c-1)$
have to be proportional  to $w.$
We have :
$${\rm Re}\ v= \lambda w +u ;\ \  {\rm Im} \ v= \gamma w,$$
with $ \lambda$ and $ \gamma\in \mathbf{R}$. If we can prove that
$\lambda =\gamma= 0$, this would imply
$v=u,$
that is $$z^a=z^b=z^c=1.$$
This in turn would mean that $z^f =1.$

Let $\rho=|z|$ be the modulus of $z$ if the rank of the matrix $A(z)$ is one then we have
 $ z^a=\lambda a+1 +i\gamma a.$
This means that  $\rho^{2a}=(\lambda a+1)^2+\gamma^2a^2$
and the same computation holds for $b$ and $c.$
This allows us to say that the following function

\begin{eqnarray*}
f:\bR & \lra & \bR,\\
 x& \mapsto & \rho^{2x} - ((\lambda x+1)^2+\gamma^2 x^2),
\end{eqnarray*}

has 4  distinct  zeroes at $a,$ $b,$ $c,$ and $0.$ The third derivative $f'''(x)$ has then the following
expression:

$$f'''(x)=(2log \rho)^3 \rho^{2x}.$$

The function $f'''(x)$ is rational over $\bR$ and by its analytical expression it is clear it has no zeroes if
$\rho \neq 1$. If $f$ has 4 zeroes, by Rolle's theorem, $f'''$ should have at least $1$ zero. This means that
$\rho=1$ and $f'''(x)$ is identically zero. Moreover we get that

\begin{equation}\label{cross}
f(x)=x(x(\gamma^2-\lambda^2)-\lambda).
\end{equation}

Equation \ref{cross} implies that for $\lambda,\gamma \neq 0$ $f$ has at most 2 zeroes but $f$ has 4 zeroes in
$0,a,b,c$. This implies that $\lambda=\gamma=0$. The converse is clear.











\end{proof}

We notice that the previous lemma for $a=1$ was proven in \cite{kaji}.
A possible generalization is the following conjecture.

\begin{conjecture}\label{matrice}

Let $0<a_1<a_2< \dots a_n$ be postive integers. Set $f:=GCD(a_1,\dots, a_n)$, then the $n\times (n-1)$ matrix

\begin{equation} 
A_{a_1,a_2,\dots,a_n}(z)=\begin{pmatrix}
 a_1  & a_1^2  &\dots & a_1^{n-2} & z^{a_1}-1 \\
a_2  & a_2^2  &\dots & a_2^{n-2}&z^{a_2}-1 \\
\dots & \dots   &\dots &\dots &\dots  \\
a_n  & a_n^2  &\dots & a_n^{n-2}&z^{a_n}-1
\end{pmatrix}
\end{equation}

should have rank smaller than $n-1$ if and only if $z^f=1$.

\end{conjecture}

The computation of this rank reduces to finding the number of positive integer zeroes of an equation of
the following type

$$z^x=p(x),$$

where $z\in \bC$ is fixed and $p(x)$ is a polynomial.

We remark that this generalization would open the way towards more general results on, for instance
 higher osculating spaces on varieties of dimension bigger than
one. A particular case of Conjecture \ref{matrice} appears in \cite{malla}, where the case $n=4$ is treated. In the paper
\cite{malla} the calculation are performed by  \textsc{MAPLE} but we
must confess that we could not decide if this process is correct.\\

We close this section giving two definitions related to Conjecture \ref{matrice}.

\begin{definition}\label{ano}
A $n$-plet $0<a_1< \dots < a_n$ of integers satisfying the condition
$GCD(a_1,\dots, a_n)=1$ such that there exists a  $z\neq 1$ such that $A_{a_1,\dots,a_n}(z)$ has rank smaller than $n-1$ is called 
an anomalous $n$-plet.
\end{definition}

\begin{definition}\label{anoc}
A locally toric curve $\phi: C\rightarrow \bP^n$ such that the $n$-plet defining its local parametrization is anomalous is called an anomalous locally toric curve.  
\end{definition}

\subsection{Rational curves}

In this section we give a geometric interpretation of Lemma \ref{lui}.\\

Consider the affine rationally parametrized curve in $\bC^3$
$$B:=\{ v(z)=(z^a,z^b,z^c),\ z \in \bC\}.$$
The only singular point of $B$ if $a>1$  is the origin.
We have $v(1)=u=(1,1,1).$  Let us consider the tangent line $L$
to the curve at $u\in B$. We have that

$$ L:=\{u+t w \}$$

where $t \in \bC$ and $w$ is like in the preceding section. Let us suppose that $L$ has another point of intersection with the curve 
$B$ and let us call it $v(\tilde{z})$.
The secant line to $B$ passing by $u$ and $v(\tilde{z})$
is $$ (1-t)u+ tv(\tilde{z})= u+t(v(\tilde{z})-u),$$

for $t\in \bC$.
This is the line $L$ if and only if $ v(\tilde{z})-u$ and $w$ are proportional, that is, using the notation of lemma \ref{lui},
 $A(\tilde{z})$ has rank $1.$ This means that Lemma \ref{lui} assures that, since $\tilde{z} \neq 1$, for such a curve with $GCD(a,b,c)=1$ the tangent
line at $(1,1,1)$ is not a trisecant.\\

Moreover we can make the following observation. Let us consider the following morphism

\begin{eqnarray*}
\mu: \bC & \lra & \bC^3,\\
x & \mapsto & (x^a,x^b,x^c).
\end{eqnarray*}

Performing for  $\lambda\in \bC ,$   $\lambda\neq 0$  the  affine transformation

\begin{eqnarray*}
\tau:\bC & \lra & \bC\\
z & \mapsto & \lambda^{-1} z
\end{eqnarray*}

on the affine line and

\begin{eqnarray*}
\delta : \bC^3 & \lra & \bC^3\\
(z_1,z_2,z_3) & \mapsto & ( \lambda^{-a} z_1, \lambda^{-b} z_2, \lambda^{-c} z_3)
\end{eqnarray*}

 on $\bC^3$, we remark that the following diagram commutes.

$$\begin{array}{ccc}
\bC & \stackrel{\mu}{\lra} & \bC^3 \\
\tau \downarrow &  & \downarrow \delta\\
\bC & \stackrel{\mu}{\lra} & \bC^3
\end{array}$$

This implies that the tangent line at any smooth point $p\in B$ does not intersect $B$ outside $p$.\\

We remark that the projective curve $\overline{B}\subset \bP^3$ such that $B$ is its restriction to an affine open set
is the image of the following morphism

\begin{eqnarray*}
\bP^1 & \lra & \bP^3,\\
{[t:z]}  & \mapsto & [ t^c:z^at^{c-a}:z^bt^{c-b}:z^c ].
\end{eqnarray*}

This means that set-theoretically we have that $\overline{B}=B\cup {s}$, where $s=(0,0,0,1).$

\section{Kaji reduction}\label{secred}

Let $C$ be a connected smooth complete algebraic curve and let

$$\phi : C \to \bP^3$$

be a morphism birational onto its image. Assume $\phi(C)$ nondegenerate.
For any point $p$ such that $\phi(p)$ is smooth, let $t_p$ be the tangent line. For any $(p,q) $
such that $\phi(p)\neq \phi(q)$ we let $s_{p,q}$ be the the secant line joining the two points.

Moreover let $V$ be the vector space $H^0(\bP^3, \cO_{\bP^3}(1))$ and $\cP^1:=\cP^1(\cO_C(1))$ the bundle
of principal parts of $\cO_C(1)$ of first order, where $\cO_C(1)=\phi^*\cO_{\bP^3}(1)$ \cite{kl}.  We recall that we have
the following commutative diagram.

\begin{equation}\label{car}
\begin{array}{ccccccccc}

0 & \ra & \phi^* \Omega^1_{\bP^3} \otimes \cO_C(1) & \ra & \cO_{\bP^3}\otimes V & \ra & \cO_C(1) & \ra & 0\\
  &     & \downarrow                               &     & \downarrow           &     & \parallel &     &  \\
0 & \ra & \Omega^1_C \otimes \cO_C(1) & \ra & \cP^1 & \ra & \cO_C(1) & \ra & 0\\
\end{array}
\end{equation}

Moreover, the bundles of principal parts of $\cO_C(1)$ of higher order realize exact sequences of the following type.

\begin{equation}\label{pp}
0 \lra Sym^m \Omega^1_C \otimes \cO_C(1) \lra \cP^m \lra \cP^{m-1} \lra 0.
\end{equation}

We recall that $\cP^1$ is not generated by sections $s \in V$ exactly over the points
of $C$ where $rk(d\phi)=0$. Let us denote $S(\phi)$ the locus where the differential has rank zero.
Let us consider now the evaluation map

$$ev:\cO_C \otimes V \lra \cP^1.$$

As we have just remarked, this map is not surjective. Let us then consider the image sheaf $E:=Im(ev)$: we have
that $E\subset \cP^1$ as a subsheaf but $E$ is also a globally generated rank 2 vector bundle. Let now
$\cO_{S(\phi)}$ be the skyscraper sheaf supported on the ramification locus of $\phi$, then we have the
following exact sequence.

$$ 0 \lra  E \lra \cP^1 \lra \cO_{S(\phi)} \lra 0$$

Furthermore we have the following commutative diagram

\begin{equation}\label{diagra}
\begin{array}{ccccccccc}
& & 0 & & 0 & & & &\\
& & \downarrow & & \downarrow & & & &\\
0 & \ra & M & \ra & E & \ra & \cO_C(1) & \ra & 0 \\
 &  & \downarrow & & \downarrow & & \parallel & & \\
0 & \ra & \Omega^1_C \otimes \cO_C(1) & \ra & \cP^1 & \ra & \cO_C(1) & \ra & 0\\
 &      & \downarrow & & \downarrow & &  & & \\
  &     & \cO_{S(\phi)} & & \cO_{S(\phi)} & & & & \\
    &     & \downarrow & & \downarrow & & & & \\
    &     & 0 & & 0 & & & & \\
\end{array}
\end{equation}

for some line bundle $M$. Let $\bP E$ be the projective bundle associated to $E$. Basically we have that $\bP E \subset
C \times \bP^3$ and the fiber $\bP E_x\subset \bP^3$ over a smooth point $x \in C$ is the tangent line to $\phi(C)$ at $\phi(x)$.

\begin{lemma}\label{pliz}
We have a map

$$\nu: \bP E \lra \bP^3.$$
\end{lemma}
\begin{proof}
We recall from diagram \ref{car} that we have a surjection

$$\cO_{\bP^3}\otimes V \twoheadrightarrow \cP^1$$

that gives us four global sections of $\cP^1$. Let $\pi: \bP E \ra C $ the natural projection and let us pull-back
diagram \ref{diagra} via $\pi$. Then we recall that $\pi^*E$ is a globally generated subsheaf of $\pi^*\cP^1$.
Let $S$ be the tautological vector bundle over $\bP E$, we are now in the following situation over $\bP E$.

$$\begin{array}{ccccccccc}
&&&& \cO_{\bP E}^{\oplus 4} &&&&\\
&&&& \downarrow &&&&\\
0 & \ra & S & \ra & \pi^* E & \ra & \cO_{\bP E^*}(1) & \ra 0\\
\end{array}$$

This implies that the four sections
of $\cO_{\bP E}^{\oplus 4}$ surject to $\cO_{\bP E^*}(1)$. Since $\pi^* E$ is globally generated, 
also $\cO_{\bP E^*}(1)$ is
globally generated, thus giving us the map we were looking for.

\end{proof}

Let us consider the map we have just defined. Let us call $Y\subset \bP^3$ the image of $\nu$, set theoretically $Y$ 
is the
closure of the union of tangent 
lines of $C$. This surface is known in the literature as the \textit{tangential surface} of
$C$.




\begin{lemma}\label{ati}
Let G an algebraic complex curve and $\cL$ a line bundle on G. Let $\cP^1(\cL)$ be the bundle of principal 
parts of $\cL$ of first order, and let
$e\in H^1(G,\Omega^1_G)$ be the extension class defined by the natural exact sequence

$$ 0 \lra \Omega^1_G \otimes \cL \lra \cP^1(\cL) \lra \cL \lra 0.$$

Then we have $c_1(\cL)=-e$. In particular,
%
$\cP^1(\cL)$ is the unique non-trivial extension of $\cL$ by $\Omega^1_G\otimes \cL$.

%
%

\end{lemma}

\begin{proof}
See \cite{aty}.
\end{proof}

Let us now consider diagram \ref{diagra}. Twisting by $\cO_C(-1)$ and passing to cohomology we get a map

$$\alpha: H^1(C,M(-1)) \lra  H^1(C,\Omega^1_C)$$

We remark that both spaces can be interpreted as parameter spaces for extension classes.

\begin{proposition}\label{pv}
Let

$$\alpha: H^1(C,M(-1)) \lra  H^1(C,\Omega^1_C)$$

be the natural map between extension classes and let a be the extension class of the first row of diagram \ref{diagra}. Then $\alpha(a)= e$.
\end{proposition}

\begin{proof}
Let us twist by $\cO_C(-1)$ diagram \ref{diagra}. We have a new diagram like the following.

$$\begin{array}{ccccccccc}\label{raf}
& & 0 & & 0 & & & &\\
& & \downarrow & & \downarrow & & & &\\
0 & \ra & M(-1) & \ra & E(-1) & \ra & \cO_C & \ra & 0 \\
 &  & \downarrow & & \downarrow & & \parallel & & \\
0 & \ra & \Omega^1_C & \ra & \cP^1(-1) & \ra & \cO_C & \ra & 0\\
 &      & \downarrow & & \downarrow & &  & & \\
  &     & \cO_{S(\phi)} & & \cO_{S(\phi)} & & & & \\
    &     & \downarrow & & \downarrow & & & & \\
    &     & 0 & & 0 & & & & \\
\end{array}$$

We remark that $M(-1) \cong \Omega^1_C(-S(\phi)$. If we take cohomology
we find the following commutative square.

$$\begin{array}{ccc}
H^0(C,\cO_C) & \lra &  H^1(C,M(-1))   \\
 \parallel & & \alpha \downarrow  \\
H^0(C,\cO_C) & \lra &  H^1(C,\Omega^1_C) \\
\end{array}$$

This implies directly the statement.

\end{proof}

Remark that Proposition $\ref{pv}$ implies that if $E$ were the trivial extension, then also $\cP^1$ would split.

\begin{definition} 
A projective curve $\phi:C \rightarrow \bP^n$ is tangentially degenerate if, for a general point $P\in \phi(C)$, there exists another point $Q\in \phi(C)$
that lies on the tangent line to $\phi(C)$ at $P$.
\end{definition}

Let us consider now the map

\begin{eqnarray*}
\rho: C \times C & \lra & C \times \bP^3, \\
(p,q) & \mapsto & (p,\phi(q)),\\
\end{eqnarray*}

and let us denote $C_0$ the image of diagonal $\Delta$ via $\rho$.

Moreover, let $\pi: \bP E \lra C$
be the natural projection, we define a map

$$\sigma: \bP E \lra C \times \bP^3$$

as the product map $\pi \times \nu$. Let us denote by $A$ the section of $\bP E$ corresponding to the surjection 
$E \rightarrow \cO_C(1)$. Remark that we have $\sigma^{-1}(C_0)=A$. Then let
$D\subset C\times C$ be the following subset:

\begin{equation}\label{fine}
D:=\{(h,s)\in C \times C \setminus C_0: t_{\phi(h)}=s_{\phi(h),\phi(s)}\}.
\end{equation}

If a curve $C$ is tangentially degenerate then $D$ has positive dimension. Moreover in this case
$\sigma^{-1}(\rho(D)) \subset \bP E$ is an effective divisor in the projective bundle.
Now we need to recall the following Proposition.

\begin{proposition}\label{gno}
Let C be a curve, let $\cE$ be a rank 2 vector bundle over $C$, and let
$\cL$ be a line bundle over $C$. Let $N$ be a section of the ruled surface $\bP \cE \ra C$ corresponding to a surjection
$\cE \ra \cL$. If there exists an effective divisor $F$ on $\bP \cE$ such that $N$ and $F$ are disjoint, then the surjection
$\cE \ra \cL$ splits
\end{proposition}

\begin{proof}
See \cite{kaji}.
\end{proof}

Furthermore, by combining Lemma \ref{ati} with Proposition \ref{gno}, 
we get the following Corollary.

\begin{corollary}\label{mdp}
Let $\phi:C \rightarrow \bP^3$ be a projective space curve. If C is tangentially degenerate, then the section $A$ of $\bP E \ra C$ intersects
every component of the curve $\sigma^{-1}(\rho(D)) \subset \bP E$.
\end{corollary}

The key point here is that if a component does not intersect $A$ then the vector bundle
$E$ splits, which is not the case. 

\begin{proposition}\label{ach}
Let $C$ be a non degenerate complex projective algebraic curve that is locally toric and 
let $\phi:C \hookrightarrow \bP^3$  be
the morphism in question. Then $C$ is not tangentially degenerate.
\end{proposition}

\begin{proof}

We suppose by contradiction that $dim(D)>0$. Let $\overline{D}$ be the closure of $D$ in $C\times C.$ We assume that $\overline{D}$ does contain an irreducible
complete curve $X.$ Let $\tilde{X}$ be the normalization of $X$ and $\psi : \tilde{X}\to C\times C$ be the
induced map. By Corollary \ref{mdp} $\sigma^{-1}(\rho(D))\cap A\neq\emptyset$. This is
equivalent
to saying that $D\cap \Delta \neq\emptyset$, hence there exists $p\in \til$ and $q\in C$ such that

$$\psi(p)=(\psi_1(p),\psi_2(p))=(q,q).$$

We take now two local analytical parameters: respectively, $t$ on an open set containing $q\in C$ and  $x$ on an open neighbourhood
of $p\in \til$. We take $t$ and $x$ such that:

\begin{enumerate}
\item [1)] $t(q)= x(p)=0$;
\item [2)] $ \psi(x)=(\psi_1(x),\psi_2(x))=(x^r, x^s(k+xg(x)),$ for a scalar $k \neq 0$ and
a function $g$ regular in a neighourhood of zero.
\end{enumerate}

We can also make a choice of the coordinates of $\bP^3$ and of the open set where the
parameter is $t$ such that we have $\phi(0)=(1,0,0,0)$ and such that in affine coordinates we have:

$$ \phi(t)=(t^a, t^b+ t^{b+1}\alpha(t), t^c+ t^{c+1}\beta(t)) $$

where $0<a<b<c$, and $\alpha$ and  $\beta$ are functions regular at zero.\\

Composing $\phi$ with $ \psi_1$ and $\psi_2$, the two components of $\psi$, we have the following parametrized curves:

\begin{eqnarray*}
\gamma (x) & = & (x^{ar}, x^{br}+x^{(b+1)r}\alpha(x^r), x^{cr}+x^{(c+1)r}\alpha(x^r)),\\
\eta (x) & = & (k^ax^{as}+higher\ terms),k^bx^{bs}+higher\  terms , k^cx^{cs}+higher\ terms),\\
\end{eqnarray*}

where $higher\ terms$ stands for the terms of higher degree.\\

Let  $\gamma' (x)$ be the first derivative of $\gamma$, then we have:

$$\gamma' (x)= (ar x^{ar-1}, brx^{br-1}+higher\ terms, crx^{cr-1}+higher\ terms),$$

since $\gamma'(x)$ gives the tangent line to $\phi(C)$ at $\phi(\psi_1(x)).$ By definition of $D$  we have:

\begin{equation} \label{setg}
D=\{(y,w)\in C\times C\setminus \Delta:(\phi(y)-\phi(w))\wedge \phi(y)'\}=0
\end{equation}

We have denoted by $\wedge$ the exterior product in $\bC^3.$ More precisely, the sub-scheme of $C\times C$
defined by the equation

\begin{equation}\label{wo}
(\phi(y)-\phi(w))\wedge \phi(y)'
\end{equation}

is the union $D\cup \Delta$, because the diagonal locus is contained in its zero locus, and scheme-theoretically with
multiplicity 2. In fact, the Taylor development of for instance $\phi(w)$ in a neighbourhood of the
diagonal gives

\begin{equation}\label{sturia}
\phi(w)=\phi(y) + (w-y) \phi'(y) +  (w-y)^2 R(y,w),
\end{equation}

where $R(y,w)$ is a holomorphic function on $C \times C$. If we substitute equation \ref{sturia} in equation
\ref{wo} we find that $\Delta$ has indeed multiplicity 2 in the zero locus of equation \ref{wo}.


Furthermore, substituting in equation \ref{wo} the lower terms of the parametrized curves
we get that the matrix

$$
B(x)=\begin{pmatrix}
arx^{ar-1},&\   x^{ar}-k^ax^{as}\\
brx^{br-1},&\   x^{br}-k^bx^{bs}\\
crx^{cr-1},&\   x^{cr}-k^cx^{cs}\\
\end{pmatrix}
$$

must have rank $1$ for all $x$. First we show that we are forced to have $s=r$. Let us suppose in fact that
$s>r$: we get that

$$
C(x)=\begin{pmatrix}
ax^{ar-1},&\   x^{ar} \\
bx^{br-1},&\   x^{br} \\
cx^{cr-1},&\   x^{cr}  \\
\end{pmatrix}
$$

has rank one for all $x$, which is false (consider for instance the first $2\times 2$
determinant). If instead $r>s$, we have the lower terms matrix

$$
D(x)=\begin{pmatrix}
ax^{as-1},&\   k^ax^{as}\\
bx^{bs-1},&\   k^bx^{bs}\\
cx^{cs-1},&\   k^cx^{cs}\\
\end{pmatrix}
$$

which has clearly rank 2 for general $x$.\\

We have then to consider only the case $ s=r.$\ In this case, by multiplying by $x$ the first column and
dividing the first (respectively, the second and the third) line by $x^{ar}$ (respectively by $x^{br}$ and
$x^{cr}$), we obtain the following matrix.

$$
A(k)=\begin{pmatrix}
 a, &\  1- k^{a}\\
 b, &\  1- k^{b}\\
 c, &\  1- k^{c}\\
\end{pmatrix}
$$

Remark that if $k=1$ the rank drops. However this corresponds to imposing the condition $y=w$ in Equation \ref{sturia},
and this condition defines $\Delta$ with multiplicity at least 2. Moreover, if we substitute $\gamma(x)$ and $\eta(x)$ in
$R(y,w)$ we find that locally in a neighbourhood of the diagonal, $D$ is defined by the equation $R(y,w)=0$. Hence the solution
$k=1$ is irrelevant to our problem because for $k=1$ we have $R(0,0)\neq 0$. By Lemma \ref{lui} the matrix $A(k)$ has rank one only if $k^f=1$,
where $f=1=GCD(a,b,c)$, so we have a contradiction. This concludes the proof.

\end{proof}

\section{Applications and Open Problems}

In this section we would like to prove, in the same spirit of the preceding section, a result about 
the span of two general tangent lines to a locally toric curve (Proposition \ref{zeus}) in $\bP^4$. Moreover we make
a conjecture about the sections of a linear system that have zeros with given multiplicity in two 
given points of a curve (Conjecture \ref{syst}). 

\subsection{Tangent lines in $\bP^4$}

\begin{remark}\label{king}
Proposition \ref{ach} states that for the general locally toric curve $C$ in $\bP^3$ we have that, keeping the same notation of section 2, $dim(D)=0$. Remark that, if one defines $D$ in the same way for a locally toric curve in $\bP^n$, we have $dim(D)=0$ as well 
for every $n$.
\end{remark} 

\begin{proposition}\label{zeus}
Let $C \subset \bP^4$ be a locally toric curve. Let $x$ and $y$ be general points of $\phi(C)$ and $z\in \phi(C)$, then the tangent lines $t_x,t_y$ and $t_z$ span the whole
$\bP^4$.
\end{proposition}

\begin{proof}
Let us consider a standard affine subset of $\bP^4$ isomorphic to $\bC^4$ and let us assume that $x,y$ and $z$ are in this subset. 



Let us now consider the product $S^2C\times C$. We suppose by contradiction  that there exists an irreducible surface $U\subset
S^2C\times C$ defined in the following way.

$$U:=\{(x,y,z)\in S^2C\times C: t_z \subset Span(t_x,t_y)\}.$$

Remark that in our construction it is not necessary to take $C\times C\times C$, in fact, since the span of $t_x$ and $t_y$ is indipendent from the order of the two points. Now we reduce this problem to the one of Proposition $\ref{ach}$, 
thus in the second part of this proof we will omit the details already stated in the proof of Proposition $\ref{ach}$. Let us consider now the \textit{partial} diagonal $\Delta'\cong C\times C$, i.e. the set of points of the type $(p + p,q), p\neq q$
 and the intersection curve $K:=U\cap \Delta'\subset \Delta'$. Remark that, since $x=y=p$, the curve $K$ is the set $D\subset C\times C\cong \Delta'$ of section 2 and the one dimensional diagonal $\Delta''$, i.e. the set of
points of the type $(p + p,p)$, is the curve $\Delta\subset C \times C$ of section 2. The fact that $dim(K)>0$ contradicts Remark \ref{king} thus concluding the proof.

\end{proof}

\subsection{The Two Osculating Points Problem}

Let us now consider again a curve $C$ of arbitrary genus and a line bundle $L$ on $C$. Let $n$ be a positive
integer and let us set a $(n+1)$-dimensional vector sub-space $V\subset H^0(C,L)$ such that the associated linear
series is base point free. We suppose moreover that the morphism

$$\varphi_V:C \ra \bP V\cong \bP^n$$

is non degenerate and birational onto its image. We recall that, given an effective divisor $F$ on $C$, the notation $V(-F)$ indicates the
intersection $V\cap H^0(C,L(-F))$ where $L(-F)$ is the kernel

$$0 \lra L(-F) \lra L \lra L_{F} \lra 0$$

of the natural evaluation map in $F$. The theory of the Wronskian assures that for a general point $p\in
C$ we have that

\begin{equation}\label{wron}
dim(V-(n+1)p)=0.
\end{equation}

This is equivalent to saying that the sub-scheme of points of $C$ such that $dim(V-(n+1)p)\neq 0$ has dimension zero. A point
such that $dim(V-(n+1)p)\neq 0$ is commonly called a \textit{flex point} \cite{gh:pag}.
Now let us consider the same problem for a couple of general points $p,q\in C$, i.e. evaluate the dimension of
the following sub-scheme of the symmetric product $S^2 C$.

\begin{definition}

We set

$$\Gamma_{n_1,n_2}:=\{(p,q)\in S^2 C : dim(V(-n_1p-n_2q))\neq 0\}.$$

\end{definition}

We remark that this definition has sense only for curves of degree $deg(C)\geq n_1+n_2$. One expects that
$\Gamma_{n_1,n_2}$ has dimension zero for $n_1+n_2=n+2$.\\

\textbf{Example} ($n=2$)\\

In this case we have a plane curve $C$ of degree at least four and we have either $n_1=1$ and $n_2=3$ or
$n_1=n_2=2$. We see that $\Gamma_{1,3}$ has dimension zero by the finiteness of flex points \cite{gh:pag} and $\Gamma_{2,2}$
is the set of bitangents, which is finite.\\

Let us now consider the case of $n=3$. In $\bP^3$ we have the sub-schemes $\Gamma_{4,1}$ and $\Gamma_{3,2}$. The
fact that $\Gamma_{4,1}$ has dimension zero is a corollary of equation \ref{wron}. 
Let $D\subset C \times C$ be the tangentially 
degenerate locus defined in Equation \ref{fine}.
For the case of $\Gamma_{3,2}$ we have the following proposition.

\begin{proposition}\label{dimensioni}
Let $C,L$ and $V$ be as before and $dim V=4$. Let us suppose that $C$ is locally toric, then

\begin{equation}\label{dim}
dim(\Gamma_{3,2})=dim(D).
\end{equation}

\end{proposition}

\begin{proof}

Let $p_t,q_t$ be points of $C$, for any $t\in \bC$. Let us suppose that we have a family of 
sections $S(t)\in V(-3p_t-2q_t)$, $t \in \bC$. 
Remark that by developing $S(t)$ we can write

$$S(t)=S(0) + tS'(0) + higher\ terms.$$

Now $S(0), S(t)\in V$, thus we have $S'(t)\in V$ too. Moreover $S(t)\in V(-3p_t-2q_t)$ and $S'(t)\in V(-2p_t-q_t)$, hence
$S(t)\cap S'(t)$ is a trisecant line that cuts out a couple of points $(p_t,q_t)\in D$. Let us now suppose that we have a family of couples 
of points $(p_t,q_t)\in D$ 
and let us take the trisecant $l_t$ line that cuts out the divisor $2p_t+q_t$. Let $p_t^*,q_t^*\in\bP^{3*}$ be the images of $p_t$ and
$q_t$ via the dual map of $C$. The dual variety of $l_t$ is an hyperlane that cuts out $2p_t^* + 3q_t^*$ on $C^*$. By repeating 
the construction of the first part of the proof then implies of a line $l_t'\in \bP^{3*}$ that cuts out $p_t^* + 2q_t^*$ on 
$C^*$. Via dualization again we find an hyperplane in $\bP^3$ that cuts out $3p_t+2q_t$ on $C$, thus concluding the proof.

\end{proof}

\begin{corollary}
Let C,L and V as in Proposition \ref{dimensioni} and let us suppose that $\varphi_V:C \rightarrow \bP V$ is locally toric, then
$dim(\Gamma_{3,2})=0$.
\end{corollary}

Let $C$ be a curve of arbitrary genus. Let moreover $L$ be as before and $V\subset H^0(C,L)$ an $(n+1)$-dimensional sub-space
such that the associated linear series is base point free. Let $\varphi_V: C \lra \bP^n$ be the associated map. More generally we have the following conjecture.

\begin{conjecture}\label{syst}\textbf{The Two Osculating Points Problem}

Let $C$ be a smooth complete complex curve such that the morphism $\varphi_V: C \lra \bP^n$ is locally toric
and let $n_1,n_2$ be two positive integers such that $n_1+n_2=n+2$. Let us moreover suppose that $deg(C)\geq
n + 2$, then the sub-scheme $\Gamma_{n_1,n_2}$ has dimension zero.

\end{conjecture}

We remark that in the case of $n_1=n$ and $n_2=2$ a costruction similar to the one used in the 
proof of Proposition \ref{dimensioni}
reduces the problem to finding osculating $(n-2)$-dimensional linear subspaces of $\bP^n$ that 
intersect the curve in a second 
point different from the osculating one. This is the problem that we study in the next section.
Remark that, for affine rational curves with an analytical parametrization of 
type

$$t \mapsto (t^{a_1},\dots,t^{a_n}),$$

proving that the general osculating space does not meet again the curve is equivalent to proving Conjecture \ref{matrice}. 

\section{A few results in higher dimension: osculating spaces in $\bP^r$}


In this section $C$ will be a locally toric curve like in section \ref{secred}, but we will consider only the genus 0 and 1 case. 
%
%
The aim of this section is to generalize Proposition \ref{ach} to the case of $(r-2)$-dimensional osculating
spaces, i.e. to show that the general osculating $(r-2)$-plane at a point of a locally toric curve does not meet the curve in a second point. When we consider non anomalous curves (see Definition \ref{anoc}) we are able to perform the reduction and show the result for every $r$. Unluckily up till now we only know that the number of anomalous $4$-plets (that are in natural bijection with local analytical parametrizations of anomalous locally toric curves in $\bP^4$) is finite and possibly zero. However it is very likely that it will be possible to show the finiteness of the anomalous $r$-plets for every $r$ (see Remark \ref{cotte}). Let $p\in C$, the projective fibre $\bP \cP^{r-2}_p$ is the osculating $(r-2)$-plane at the point $p$. 
Let us give the following definition,

\begin{equation}\label{di}
D_{r-2}:=\{(p,q)\in C\times C\setminus  \Delta : q \in \bP \cP^{r-2}_p\}.
\end{equation}

Before giving the main theorem of this section we need some technical Lemmas.
Let now $(\alpha,\beta)\in C\times C$ and let $\Delta_{\alpha,\beta}$ the translated diagonal

$$\Delta_{\alpha,\beta}:=\{(p+\alpha,p+\beta)\in C\times C, \forall p \in C\}.$$

\begin{lemma}\label{int}
Let $C$ be a genus 1 curve. The curves $\Delta_{\alpha,\beta}$ are the only irreducible 
curves in $C\times C$ that do
not intersect the diagonal.
\end{lemma}

\begin{proof}
Let us suppose that $\Delta_{\alpha,\beta}\cap \Delta\neq\emptyset$. Hence there exist two 
points 
$p,q\in C$ such that $p +\alpha = q$ and $p+\beta=q$. This is true only if $\alpha=\beta=0$ and $p=q$.
Remark that, since $C\times C$ is an abelian variety, every divsor in $C\times C$ is NEF, thus it is ample or it has autointersection equal to zero. The curve $\Delta_{\alpha,\beta}$ is not ample and this concludes the proof.

\end{proof}

\begin{lemma}\label{como}
Let C be an elliptic curve, $L\in Pic^0(C)$ and $\varphi_V:C \rightarrow \bP V$ as before. Let $\cP^{n*}$ be the bundle of principal parts of order $n$ of $\cO_C(1)$ (see diagrams \ref{car} and \ref{pp}), then

\begin{eqnarray}
h^0(C,\cP^{n*}\otimes \cO_C(1) \otimes L)= 0 & \mathrm{if}\ L \neq \cO_C;\\
h^0(C,\cP^{n*}\otimes \cO_C(1) \otimes L)\neq 0 & \mathrm{if}\ L = \cO_C.
\end{eqnarray}

\end{lemma}

\begin{proof}
By dualizing and twisting by $\cO_C(1)\otimes L$ the second row of diagram \ref{car} we get 

\begin{equation}\label{bottle}
0 \longrightarrow L \longrightarrow  \cP^{1*} \otimes \cO_C(1)\otimes L \longrightarrow - \Omega^1_C \otimes L \longrightarrow 0.
\end{equation}

If $L \neq \cO_C$ then $h^0(C,L)=h^0(C,-\Omega^1_C \otimes L)=0$ and passing to cohomology we find that
$h^0(C,\cP^{1*}\otimes \cO_C(1) \otimes L)= 0$. Now, by dualizing and twisting by $\cO_C(1)\otimes L$ the exact sequence \ref{pp}, by induction we get that $h^0(C,\cP^{n*}\otimes \cO_C(1) \otimes L)= 0$.  If $L=\cO_C$ then $h^0(C,L)=h^1(C,L)=1$. Remark that the coboundary map 

$$H^0(C,-\Omega^1_C\otimes L)\cong H^0(C,L)  \longrightarrow  H^1(C,L)$$

is an isomorphism. This implies that $h^0(C,\cP^{1*}\otimes \cO_C(1) \otimes L)= 1$. Now, by dualizing and twisting by $\cO_C(1)\otimes L$ the exact sequence \ref{pp}, by induction we get that $h^0(C,\cP^{n*}\otimes \cO_C(1) \otimes L)>0$.
\end{proof}

\begin{theorem}\label{grosso}
Let $C\subset \bP^r=\bP V$ be as before in this section and $g(C)=0$ or $1$. Moreover, if $g(C)=1$ we suppose that the the evaluation map 

$$ev_{r-2}:V \longrightarrow \cP^{r-2}$$

is surjective. If $C$ is non anomalous then $dim(D_{r-2})=0$, i.e. the general 
osculating $(r-2)$-plane does not intersect the curve in a second point different from the osculating one.
\end{theorem}

\begin{proof}

By contradiction we suppose that $D_{r-2}$ has dimension 1. Let $\tilde{D}$ be the image in $C\times C$ of the normalization 
of one irreducible component of $D_{r-2}$. Let $\Delta$ be as usual the diagonal. If we are able
to show that $\tilde{D}\cap \Delta \neq \emptyset$ then, since we suppose that $C$ is non anomalous, we have a contradiction. In fact this would imply that, by repeating the proof of Proposition \ref{ach}, that the matrix (\ref{Amm}) has rank $r-2$ for all $z$.
This is very easy when $g(C)=0$. In fact it is well known that in this case 
$\Delta$ is ample as a divisor of $C\times C$ and thus it intersects every curve in the product. Now we come to the genus 1 case. Here we have that
$\Delta^2=2-2g=0$ and $\Delta$ is NEF. Now let us suppose that the curve $\tilde{D}$ has no point in common with the 
diagonal $\Delta$. By Lemma \ref{int} this means that there exist two points $\alpha,\beta\in C$ such that 
$\tilde{D}=\Delta_{\alpha,\beta}$.

Let us denote $\varphi_V: C \rightarrow \bP V$ the locally toric morphism.
We have a rank 2 vector bundle over $\tilde{D}$ defined in the following way:

$$\bP N_{(p,q)}:= s_{\varphi_V(p),\varphi_V(q)},$$

i.e. the projectivized fiber over every couple of points in $\tilde{D}$ is the secant line joining the two points. More precisely, let us denote $Gr(2,r+1)$ the Grassmannian of 2-spaces in $V$ and let

$$\tau: C \times C /\Delta \longrightarrow Gr(2,r+1)$$

be the secants map. We have the tautological rank 2 vector bundle $T$ on $Gr(2,5)$ and
$N$ is the restriction of $\tau^*T$ to $\tilde{D}$. Now let
us denote $\pi_1$ (respectively $\pi_2$) the projection from $\tilde{D}$ on the first factor (respectively the second). We have
a natural inclusion $\bP N \hookrightarrow  \bP \pi_1^* \cP^2$ that comes from the fact that if $(p,q)\in \tilde{D}$ then the secant  $s_{\varphi_V(p),\varphi_V(q)}$ is contained in $\bP \cP_p^2$. This is equivalent to a surjection

$$\pi_1^* \cP^2 \stackrel{g}{\longrightarrow}  N \longrightarrow 0.$$

Moreover we have two natural sections $D_1$ and $D_2$ of $\bP N$ over $\tilde{D}$. In fact the curve $D_1$ (resp. $D_2$) represents the image in $\bP N$
of the first (resp. the second) component of $\tilde{D}\subset C \times C$. Thus we have two surjections

$$\pi_1^* \cP^2  \stackrel{s_1}{\longrightarrow}   \pi_1^*\cO_C(1)  \longrightarrow 0,$$
$$\pi_1^* \cP^2  \stackrel{s_2}{\longrightarrow}   \pi_2^*\cO_C(1) \longrightarrow 0,$$

and both factorize via $g$. This implies that we have a third surjection

\begin{equation}\label{edge}
\pi_1^* \cP^2 \stackrel{s_1 \oplus s_2}{\longrightarrow} \pi_1^*\cO_C(1) \oplus \pi_2^*\cO_C(1) \longrightarrow 0.
\end{equation}

We remark that since $\tilde{D}=\Delta_{\alpha,\beta}$ we have forcely 
$deg(\pi_1)=deg(\pi_2)=1$ and
thus $deg(\pi_1^*\cO_C(1))=deg(\pi_2^*\cO_C(1))= deg(\cO_C(1))$. Furthermore we can write $\pi_2^*\cO_C(1)$ as
$\pi_1^*\cO_C(1)\otimes L$, where $L \in Jac(\tilde{D})$. Remark that $L\neq\cO_C$ since $\alpha\neq\beta$.



Now Lemma \ref{como} implies that the last arrow of the exact sequence \ref{edge} cannot be a surjection.
This in turn implies that $\tilde{D}$ cannot be a curve of type $\Delta_{\alpha,\beta}$ and thus 
the intersection $\tilde{D}\cap \Delta$ is not empty. This concludes the proof for $g(C)=1$. 

\end{proof}



\bibliographystyle{amsplain}
\bibliography{bibkaji}



Michele Bolognesi\hfill   Gian Pietro Pirola\\
Scuola Normale Superiore\hfill  Dipartimento di Matematica ''F. Casorati"\\
P.za dei Cavalieri, 7\hfill   Via Ferrata, 1\\
Pisa, 56100, ITALY\hfill   Pavia, 27100, ITALY\\
michele.bolognesi@sns.it\hfill  gianpietro.pirola@unipv.it\\

\newpage

\newcount\notenumber
\def\clearnotenumber{\notenumber=0}
\def\note{\advance\notenumber by 1
\footnote{$^{(\the\notenumber)}$}}
\font\title=cmr10 scaled 1200

\def\b{{\bf b}}
\def\A{{\bf A}}
\def\Z{{\bf Z}}
\def\R{{\bf R}}
\def\Q{{\bf Q}}
\def\C{{\bf C}}
\def\c{{\cal C}}
\def\O{{\cal O}}
\def\N{{\bf N}}
\def\Pr{{\bf P}}
\def\Gm{{\bf G}_m}
\def\x{{\bf x}}
 
\font\title=cmr10 scaled 1200

\centerline{{\title APPENDIX:
 On the rank of certain matrices}}
\bigskip

\centerline{Pietro Corvaja \& Umberto Zannier}
\bigskip

$${}$$

The aim of this note is to prove the following result, motivated by Conjecture 1.2 of Bolognesi and Pirola [B-P].
\medskip

{\bf Theorem}. {\it There are only finitely many points $(\xi,a,b,c,d)\in\C\times\Z^4$, with $\xi\in\C, |\xi|>1$, $0<a<b<c<d$, $\gcd(a,b,c,d)=1$, such that

$${\rm rank} \begin{pmatrix}
 1  & 0  & 0 & 1 \\
\xi^a  & a^2  & a & 1 \\
\xi^b & b^2   & b & 1  \\
\xi^c  & c^2 & c & 1\\
\xi^d  & d^2 & d & 1
\end{pmatrix}<4 \eqno(*)$$


Also, these points can be effectively determined.}
\bigskip

Clearly, the condition of coprimality is just a normalization condition: from a solution $(\xi,a,b,c,d)$ with $h=\gcd(a,b,c,d)$ one gets the ``primitive" solution\\
$(\xi^h, a/h,b/h,c/h,d/h)$ and conversely.\medskip

The underlying geometrical problem appeared also in papers previous  to [B-P] (see the references therein). This also leads to an analogous problem for  $(n+1)\times n$ matrices, whose natural formulation is given in [B-P, Conjecture 1.2], the present case corresponding to $n=4$. The case $n=2$ is clear, whereas for  $n=3$ Bolognesi  and Pirola show by the following  simple but ingenious argument that there are indeed no solutions at all: any solution  gives rise to equations  $\xi^m=f(m)$ for $m=0,a,b,c$, where $f(x)=\alpha x+\beta$ is a suitable  linear complex polynomial. Taking complex conjugates and multiplying yields $|\xi|^{2m}=q(m)$ for the same values of $m$, where $q=|f|^2$ is a real polynomial  of degree $\le 2$. But now Rolle's theorem, applied three times to the real function $|\xi|^{2x}-q(x)$ with four distinct zeros,  shows this is impossible. (A similar proof has been given independently also by Elkies in January 2005, as an answer to a question by Izzet Coskun, and Joseph Harris, and we thank Elkies for forwarding us his solution.) This method however does not work for $n>3$. 

Conjecture 1.2 of [B-P] states that for any $n$ there are no solutions, i.e. no  $(n+1)\times n$ matrices as in (*) with rank $<n$. Our theorem reduces the case $n=4$ to a finite computation (see the remarks at the end on how this can possibly be implemented in practice). We believe that the present method is able to yield an analogue of the above theorem for every $n$, however with a finite number of exceptions depending possibly on $n$. 
\medskip

The proof of the theorem  is partially based on the techniques appearing  in [B-M-Z].

We shall adopt the usual notation $h(\cdot)$ for the logarithmic Weil height as for instance in [B-G]. 
We shall also need the following auxiliary results:
\bigskip

{\bf Lemma 1} (Weak form of Dobrowolski). {\it For every $\epsilon>0$ there exists a computable number $c(\epsilon)>0$ such that for every algebraic number $\xi$, not zero and not root of unity,   
$$
h(\xi)>c(\epsilon)[\Q(\xi):\Q)]^{-1-\epsilon}.
$$
}
\medskip

The following result appears in [Z], Lemma 1:\smallskip

{\bf Lemma 2}. {\it Let $\xi,\gamma_1,\ldots,\gamma_h\in\overline{\Q}$, $m_0>m_1>\ldots>m_h=0$ be integers with
$$
\xi^{m_0}+\gamma_1\xi^{m_1}+\ldots+\gamma_h=0.
$$
Assume also that no subsum of the form 
$ \xi^{m_0}+\ldots+\gamma_l\xi^{m_l}$, $l<h$, vanishes.
Then there exists a computable number $B_1(h)$ such that
$$
h(\xi)\leq B_1(h){\max_i(h(\gamma_i))+1 \over m_0}.
$$
} 
\bigskip

{\it Proof of Theorem}. In the sequel $B_2,B_3,\ldots$ denote computable absolute constants. The system coming from $(*)$ in the statement consists of five quadrinomial equations in $\xi$; we do not write them explicitly,  but we only remark that: 

(i) all the involved coefficients are nonzero and have height $\ll \log d$;

(ii) the exponents for $\xi$ in the five equations are respectively the five quadruples in the set $\{0,a,b,c,d\}$.
\medskip

{\bf Step 1}. As a first step,  applying the above Lemma 2 to each of the said quadrinomial equations, one can prove the bound
$$
h(\xi)\leq B_2\cdot {\log d\over d}\eqno(1)
$$
for all solutions $\xi$ to the system $(*)$.

We sketch the argument.
This follows immediately if no relevant subsum vanishes in at least one of the three equations containing both exponents $0, d$.  
So we may assume that for each of these three equations we have some ``initial subsum" vanishing. Note that since no coefficient vanishes, such subsum  must be binomial. 

Taking the equation with exponents $0,a,b,d$, we obtain (from the binomial involving $b,d$) an estimate similar to $(1)$ with $d-b$ in place of $d$ in the denominator. Then, arguing with the exponents $0,b,c,d$,  we obtain (from the binomial involving $0,b$) a similar estimate with $b$ in place of $d$ in the denominator. 

Hence in any case we have a bound similar to (1), with $\max (b,d-b)\geq d/2$ in place of $d$, concluding the proof of $(1)$.
\medskip

Now, given four coprime integers $0<a<b<c<d$, let us consider the variable vector $\x=(x_a,x_b, x_c,x_d)$ and  define the linear variety $L$ by the condition

$${\rm rank} \begin{pmatrix}
 1 & 0 & 0 & 1 \\
x_a & a^2 & a & 1 \\
x_b & b^2 & b & 1  \\
x_c & c^2 & c & 1 \\
x_d & d^2 & d & 1
\end{pmatrix}<4$$


Clearly the variety $L$ is a plane in ${\bf A}^4$. With respect to the present coordinates it is
``non degenerate" in the sense that it is not contained in any affine hyperplane defined by a linear equation with three or less terms. This follows from the Van der Monde structure of the matrix in question. In fact, the plane $L$ can be parametrized as follows:   the rank is $< 4$ if and only if the first column is a linear combination of the last three. Namely,   ${\bf x}=(x_a,x_b,x_c,x_d)\in L$ if and only if there exist $(x,y)\in\C^2$ such that 
$$
x_a=a^2x+ay+1, x_b=b^2x+by+1, x_c=c^2x+cy+1, x_d=d^2x+dy+1.\eqno(2)
$$
These equations also yield a parametrization of the plane $L$ by the parameters $x,y$. 
In particular,  no two  of the coordinates $x_a,\ldots,x_d$ satisfy a linear equation identically in $L$.\medskip

We are interested in the points $P_\xi=(\xi^a,\xi^b,\xi^c,\xi^d)$ in $L$ for a complex number $\xi$ with $|\xi|>1$.

\medskip

In order to apply Lemma 1, we need an upper bound for the degree of $\xi$, which shall be obtained via geometry of numbers. Let $\|\cdot \|$ denote any norm on $\R^4$, e.g. the euclidean length;  let $\Lambda\subset  \Z^4$ be the lattice of vectors  orthogonal to $(a,b,c,d)$. Minkowski's second Theorem [B-G, Thm. C.2.11] provides the  existence of a basis   $\lambda_1,\lambda_2,\lambda_3$  for $\Lambda$  satisfying
$$
\|\lambda_1\|\leq\|\lambda_2\|\leq \|\lambda_3\|; \qquad \|\lambda_1\|\cdot\|\lambda_2\|\leq B_3  d^{2/3},
$$
with a constant $B_3$ depending only on the chosen norm. 
We remark that  a careful choice of the norm may produce values of $B_3$ more convenient for actual computations, as we shall observe at the end.\medskip

Let now $G\subset\Gm^4$ be the  (two dimensional) torus defined by
$$
{\bf x}^{\lambda_1}={\bf x}^{\lambda_2}=1.
$$
This torus corresponds to the primitive lattice generated by $\lambda_1,\lambda_2$, in the sense that a point ${\bf x}$ is in $G$ if and only if ${\bf x}^\lambda=1$ for every $\lambda$ in the said lattice.

Let us suppose first that $L\cap G$ is finite, as expected by dimensional considerations. In this case, by  B\'ezout's theorem, $\sharp(L\cap G)\leq B_4 d^{2/3}$; since   both varieties $L$ and $G$ are defined over $\Q$,  the degree (over $\Q$) of the point $(\xi^a,\xi^b,\xi^c,\xi^d)$ is $\leq B_4 d^{2/3}$. The coprimality of $a,b,c,d$ implies the same bound for the degree of $\xi$:
$$
[\Q(\xi):\Q]\leq B_4 d^{2/3}.
$$ 
Now Lemma  1 with any $\epsilon<1/2 $ provides an explicit upper bound for $d$, concluding the proof in this case.   
\medskip

Even if $\dim(L\cap G)\geq 1$  the above method  applies, provided  $P_\xi$ is an   isolated component of $L\cap G$.
\smallskip

Let then assume that $P_\xi$ lies in a component $\c$ of $L\cap G$ of positive dimension.
This situation of anomalous dimension of the intersection of subtori with algebraic subvarities of  $\Gm^N$ is the object of [B-M-Z]. We follow some of the arguments therein.\medskip

 {\bf Step 2}. {\it We prove that $\c$ is a conic or a line in $L$. }

The non degeneracy of $L$ implies that $\c$ has dimension $\leq 1$: otherwise, $L$ would be equal to $G$, so in particular $x_a,x_b,x_c$ would be multiplicatively dependent on $L$, but the above representation $(2)$ excludes this, by the non-proportionality of the linear forms which appear. 
So we can suppose that $\c$ is a curve.

Since the equations of the plane $L$ and the equations $\x^{\lambda_1}=1, \x^{\lambda_2}=1$ hold identically on $\c$, the corresponding jacobian  determinant  vanishes on $\c$.

Writing $\lambda_{i}=(l_{i,a},l_{i,b},l_{i,c},l_{i,d})$ for $i=1,2$, the equations for $G$
read: $x_a^{\lambda_{1,a}}\cdots x_d^{\lambda_{1,d}}=1; x_a^{\lambda_{2,a}}\cdots x_d^{\lambda_{2,d}}=1$. Equations for $L$ can be written as $x_c=\beta x_a+\gamma x_b -\delta$, $x_d=\beta^\prime x_a+\gamma^\prime x_b-\delta^\prime $.
 We arrive at the equation

$$det \begin{pmatrix}
\beta&\gamma&-1&0 \\
\beta^\prime&\gamma^\prime&0&-1 \\
\frac{\lambda_{1,a}}{x_a}&\frac{\lambda_{1,b}}{x_b}&\frac{\lambda_{1,c}}{ x_c}&\frac{\lambda_{1,d}}{x_d}\\
\frac{\lambda_{2,a}}{x_a}&\frac{\lambda_{2,b}}{x_b}&\frac{\lambda_{2,c}}{ x_c}&\frac{\lambda_{2,d}}{x_d}
\end{pmatrix} =0$$


Due to the nondegeneracy of $L$, i.e. equations (2), no $2\times2$ minor of the first two rows vanishes.
So, by the linear independence of $\lambda_1,\lambda_2$, this determinant does not vanish identically as a function of four independent variables. Actually, substituting $x_a,x_b,x_c,x_d$ by their representation $(2)$ valid on $L$, the determinant does not vanish identically on $L$, since the linear polynomials in $(2)$ are pairwise non proportional.
After clearing the denominators and using the parametrization (2), this becomes a non trivial quadratic equation in $x,y$. We then obtain the statement of the Step 2.
\bigskip

{\bf Step 3}. {\it Conclusion}. If $\c$ is a line, then each coordinate $x_a,\ldots,x_d$ is a linear function of a parameter $t$; by the parametrization $(2)$, at most one of them can be constant on $\c$. Take any equation ${\bf x}^\lambda=1$ valid on $G$, containing a coordinate non constant  on $\c$. Such a coordinate will vanish on a point $P\in\c$, so  the equation must contain another coordinate vanishing on $\c$, 
 necessarily proportional on $\c$ to the first one. 
On the other hand, we cannot have three pairwise proportional coordinates on $\c$, because otherwise three of the linear forms in $(2)$ would have a common zero. In view of the previous conclusion, and since $G$ has codimension two, the four coordinates must be proportional in pairs, in particular non constant. Now, the ratio between two multiplicatively dependent and proportional linear polynomials is necessarily a root of unity.
But then $P_\xi$ cannot lie on $\c$, otherwise $\xi$ itself would be a root of unity.\smallskip

Suppose now $\c$ is a conic.
We argue as follows:
 first note that  in this case no coordinate function can be constant on $\c$, nor can two coordinates  be proportional on $\c$, otherwise $\c$ would be a line.
Then we distinguish two cases: ({\it P}) $\c$ is a parabola and ({\it H}) $\c$ is a hyperbola.\smallskip

In case ({\it P}), all the four coordinates (viewed as regular functions on $\c$) have a pole at the unique point at infinity of $\c$; suppose now that two coordinates, say $x_a,x_b$, have a simple pole at infinity; then the same must hold also for $x_c,x_d$,   since from the representation (2) it follows that $x_c,x_d$ are linear combinations of $x_a,x_b$. 
But this is impossible, since in this case  $\c$ would be a line. Then it remains to consider the case when at most one coordinate has a simple pole, so the other three, say $x_a,x_b,x_c$ have each a double pole at infinity.
Let us write  a nontrivial multiplicative dependence relation, which, after possibly permuting the indices, takes the form $x_a^mx_b^n=x_c^l$, with $m,n,l$ non negative and not all zero (such a relation  exists since $G$ has dimensione two). Looking at the poles we get $l=m+n$; looking at the zeros, we obtain that necessarily two out of $x_a,x_b, x_c$ have the same zeros. Then the corresponding ratio would be constant, contrary to the assumption that $\c$ is not a line.

Let us now consider case ({\it H}), when $\c$ is a hyperbola, hence has two points $Q_1\neq Q_2$ at infinity. We distinguish three subcases:   ({\it Ha}) there exist two coordinates, multiplicatively independent on $\c$,  having  poles both at $Q_1$ and at $Q_2$. 
 ({\it Hb}) any two coordinates having poles both at $Q_1$ and  $Q_2$ are multiplicatively dependent on $\c$; ({\it Hc}) there is at most one coordinate having poles both in $Q_1$ and in $Q_2$.

In case ({\it Ha}), we can suppose $x_a,x_b$ are multiplicatively independent on $\c$ and have poles both at $Q_1$ and at $Q_2$. Since the group $G$ has dimension $2$, there exist positive powers of both $x_c$ and $x_d$  which equal a product of powers of $x_a,x_b$ on $G$, hence on $\c$. From a relation $x_c^e=x_a^mx_b^n$, with $e,m,n\in\Z,\, e>0$ it follows that 
$m+n=e$ and that $x_c$ too has both poles, other wise $x_c$ would have no pole, so it would be constant which we have excluded. Then, after interchanging if necessary the roles of $x_a,x_b,x_c$, we can suppose $m>0, n\leq 0$. If $n=0$, then $x_a/x_c$ would be constant, which is impossible as observed. If $n<0$, then each zero of $x_b$ (on the conic $\c$) would be a zero of $x_a$ (otherwise $x_c$ would have another pole); now, if $x_b$ had two simple zeros, $x_a$ should have exactly the same zeros with the same multiplicities, so $x_a/x_b$ would be constant, which is excluded. Then $x_b$ would have a double zero, say $P_1$, $x_c$ another double zero, say $P_2$, and $x_a$ has zeros in $P_1,P_2$. Repeating the same argument with $x_d$, we see that the zero set of $x_d$ is also contained in $\{P_1,P_2\}$. Since there are only three possibilities for the zero divisor of the four functions $x_a,x_b,x_c,x_d$, namely $P_1+P_2, 2P_1, 2P_2$, some ratio must be constant.

In case {\it (Hb)}, suppose $x_a,x_b$ have (simple) poles in $Q_1,Q_2$ and are multiplicatively dependent. But then from a  relation $x_a^mx_b^n=1$ with $(m,n)\in\Z^2\setminus\{0\}$ it follows looking at the poles that $n=-m$, so $x_a/x_b$ is constant, which would imply that $\c$ is a line.

Finally, in case {\it (Hc)}, two coordinates would have poles just in one point $P_1$; but  we have already remarked that by (2) any two coordinates generate the others as linear combinations, so all coordinates would have just one points at infinity, the same for all, contrary to the assumption that $\c$ is a hyperbola.
\medskip

This concludes the sketch of the proof. We remark that Step 3 does not follow [B-M-Z]. An alternative but more involved argument would follow [B-M-Z] more closely: using the fact that $\deg(\c)$ is bounded, we would first prove that $\lambda_1$ has bounded length. This would allow a finite number of parametrizations of $a,b,c,d$ in terms of three unknowns which would be the first step for an induction. Perhaps in case of higher dimension this last method will be necessary.
\bigskip

\noindent {\bf \S 2 Quantification}. 
\medskip

Fact 1. In place of Dobrowolski Theorem (Lemma 1), we can use Smyth result: \medskip

{\bf Lemma 3}. {\it Let $\xi$ be a non reciprocal algebraic number. Then 
$$
h(\xi)\geq \frac{\log \theta_0}{[\Q(\xi):\Q]}\geq \frac{0.28}{[\Q(\xi):\Q]}
$$
where $\theta_0$ is the real root of $z^3=z+1$.
}
\medskip

In fact, we prove as follows that we can disregard reciprocal numbers $\xi$.
If $(*)$ is verified there exists a non zero polynomial $f(X)\in\bar{\Q}[X]$ of degree $\leq 2$ such that
$f(m)=\xi^m$ for $m=0,a,b,c,d$. If $\xi$ is reciprocal, there exists a Galois automorphism $\sigma$ such that $\xi^\sigma=\xi^{-1}$. Then applying $\sigma$ and multiplying we obtain that the polynomial $f(X)f^{\sigma}(X)-1$ has five distinct zeros which may happen only if $f$ is constant, which is a case easily dealt with.
\medskip

A more explicit version of Lemma 2 will lead to explicit quantitative bounds:\medskip

{\bf Lemma 4}. {\it Let $\xi,\gamma_1,\ldots,\gamma_h\in\overline{\Q}$, $m_0>m_1>\ldots>m_h=0$ be integers with
$$
\xi^{m_0}+\gamma_1\xi^{m_1}+\ldots+\gamma_h=0.
$$
Let $l<h$ be such that the subsum
$$
\xi^{m_0}+\gamma_1\xi^{m_1}+\ldots+\gamma_l \xi^{m_l} \neq 0
$$
Then 
$$
(m_l-m_{l+1})h(\xi)\leq h(1:\gamma_1:\ldots:\gamma_h) +\log\max\{l+1,h-l\}.
$$
}
\medskip

{\it Proof}. One applies the product formula $\prod_\nu|\phi|_\nu=1$to $\phi:=\xi^{m_0}+\gamma_1\xi^{m_1}+\ldots+\gamma_l \xi^{m_l}$. 
We estimate the various factors as follows: if $|\xi|_\nu<1$ then
$$
|\phi|_\nu<  |\xi|_\nu^{m_l}\sup(1,|\gamma_1|_\nu,\ldots,|\gamma_l|_\nu)\cdot\sup(1, |l+1|_\nu).
$$
If $|\xi|_\nu\geq 1$, one uses the equation $\phi=-(\gamma_{l+1}\xi^m_{l+1}+\ldots+\gamma_h)$ to obtain
$$
|\phi|_\nu \leq |\xi|_\nu^{m_{+1}}\sup(1,|\gamma_{l+1}|_\nu,\ldots,|\gamma_h|_\nu)\cdot\sup(1, |h-l|_\nu).
$$
Taking the product over all $\nu$ one obtains the estimate of the Lemma.

\bigskip

The three equations coming from $ (*)$ and involving the term $\xi^d$ are of the shape
$$
\xi^d rs(s- r)\pm \xi^s dr(d-r)\pm \xi^r ds(d-s) \pm (d-r)(d-s)(s-r)=0,
$$
where $r<s$ are two integers in $\{a,b,c\}$.

If, as in Step 1, in at least one of the equations no initial subsum vanishes,
  we can apply Lemma 4 for $l=0,1,2$  to find

$$
d h(\xi)\leq 9 \log d -\log {8^3\over 18}\leq 9\log d - 3.347.
$$

This inequality holds {\it a fortiori} in the other cases (by following Step 1). 
\medskip

To quantify the step involving Minkowski's Theorem,  
one can check that 
$$
\|\lambda_1\| \cdot\|\lambda_2\|\leq {96 d}^{2/3}.
$$
Here the involved norm $\|\cdot\|$ is the sum of the absolute values of the coordinates; this choice of the norm is motivated by the fact the degree of the closure of the variety $x^{\lambda}=1$ is $\leq \|\lambda\|$.  So the bound $96 d^{2/3}$ holds for the degree of $\xi$.

Finally we obtain
$$
\frac{0.28\cdot  d^{1/3}}{96^{2/3}}\leq 9\log d-\log \frac{8^3}{18}.
$$
This bound is explicit, for instance it gives $d<10^{13}$, but this is too  large for explicit computations of all the relevant cases.

To obtain more realistic estimates one could proceed as follows:\medskip

Let $a_0\xi^t+\ldots +a_t=0$ be a minimal equation for $\xi$ over $\Z$ (so in particular $t=[\Q(\xi):\Q]$). Then trivially one has $h(\xi)\ge (\log\max(|a_0|,|a_t|))/t$. This largely improves on Lemma 3 unless $a_0,a_t$ are both small, i.e. $\xi$ is `almost' a unit. (Note that already if $\max (|a_0|,|a_t|)\ge 2$ we replace $0.28$ by $\log 2>0.69$.) 

If this is the case however we can gain on refining Lemma 4. For instance note that if $\xi$ is actually a unit using the above equations we find for each conjugate $|\xi^\sigma|^{d-s}\le\approx {d(d-r)\over s(s-r)}$. For a `general' choice of $r,s$ this ratio on the right  will be absolutely bounded, providing a further saving.

Finally, the right norm to consider in applying Minkowski comes from an estimation for the degree of the closure  in $\Pr_4$ of a hypersurface defined by $x_1^{m_1}\cdots x_4^{m_4}=1$, where $m_i\in\Z$. The degree equals the maximum between the sum of the positive, or $-$ negative, entries of $(m_1,\ldots ,m_4)$. This quantity defines actually a norm, and yields better constants than the above  choice.\medskip

To perform this program however leads to several computing verifications and distinction into cases, which we have by now not performed, hoping that some expert of such matters can be interested enough to do this job. In this respect we remark that for the geometric application in the paper, it is substantial to pass from finitely many solutions to no solution  at all.

\bigskip

\noindent{\bf References}
\medskip

[B-G] E. Bombieri, W. Gubler, Heights in Diophantine Geometry, Cambridge U. Press 2006.\smallskip

[B-M-Z] E. Bombieri, D. Masser, U. Zannier, Anomalous subvarieties, to appear in {\it Int. Math. Res. Notices}.\smallskip

[B-P] M. Bolognesi, P. Pirola, Osculating spaces and diophantine equations, preprint 2007.\smallskip

[Z] U. Zannier, Appendix to the book by A. Schinzel: Polynomials with special regard to reducibility, Encyclopedia of mathematics and its applications, Cambridge U. Press 2000

\vfill

$${}$$

Pietro Corvaja\hfill Umberto Zannier

Dipartimento di Matematica e Informatica\hfill  Scuola Normale Superiore

Via delle Scienze, 206\hfill Piazza dei Cavalieri, 7

33100 - Udine (ITALY)\hfill  56100 Pisa (ITALY)

corvaja@dimi.uniud.it\hfill u.zannier@sns.it


\end{document}